\documentclass[12pt,leqno]{article}

\title{\bf\normalsize DEFINABLE DAVIES' THEOREM}
\author{\normalsize Asger T\" ornquist\footnote{A. T\" ornquist was supported in part by the Danish
Natural Sciences Research Council post-doctoral grant no.
272-06-0211}, William Weiss}

\usepackage{amsmath}
\usepackage{amstext}
\usepackage{amssymb}
\usepackage{amsthm}
\usepackage{enumerate}
\usepackage{mathrsfs}

\DeclareMathOperator{\lh}{lh}

\DeclareMathOperator{\supp}{supp}\DeclareMathOperator{\dom}{dom}

\DeclareMathOperator{\IS}{IS} \DeclareMathOperator{\Succ}{succ}
\DeclareMathOperator{\lex}{lex}

\newcommand{\restrict}{\upharpoonright}

\DeclareMathOperator{\rank}{rank}

\def\R{{\mathbb R}}

\def\N{{\mathbb N}}

\def\({{\normalfont (}}
\def\){{\normalfont )}}

\newcounter{mypar}
\newcounter{thmcounter}

\newcommand{\mysec}[1]{\setcounter{thmcounter}{0}\addtocounter{mypar}{1}\section*{\begin{center}\normalsize{\sc \S \arabic{mypar}.
#1}\end{center}}}

\newcommand{\mythm}[2]{\addtocounter{thmcounter}{1}\subparagraph{{\sc \arabic{mypar}.\arabic{thmcounter}. #1}}{\it
#2}}

\newcommand{\thm}[2]{\subparagraph{{\sc#1\ }}{\it
#2}}

\newcommand{\remark}[1]{\subparagraph{\it #1}}

\newcommand{\myheadpar}[1]{\addtocounter{thmcounter}{1}\subparagraph{\rm \arabic{mypar}.\arabic{thmcounter}. {\it #1}}}

\newcommand{\claim}[1]{\subparagraph{{\sc #1}}}

\begin{document}

\textwidth=30cc \baselineskip=16pt

\maketitle

\begin{abstract}
We prove the following descriptive set-theoretic analogue of a
Theorem of R.O. Davies: Every $\Sigma^1_2$ function
$f:\R\times\R\to\R$ can be represented as a sum of rectangular
$\Sigma^1_2$ functions if and only if all reals are constructible.
\end{abstract}

\mysec{Introduction}

In \cite{davies}, R. O. Davies proved that the continuum hypothesis,
CH, is equivalent to the statement that every function
$f:\R\times\R\to \R$ can be represented as a sum of ``rectangular''
functions as follows: There are $g_n,h_n:\R\to\R$, $n\in\omega$,
such that
$$
f(x,y)=\sum_{n=0}^\infty g_n(x)h_n(y),
$$
where at each $(x,y)\in\R^2$ there are at most finitely many
non-zero terms in the above sum. We call such a representation a
\emph{Davies representation} of $f$. Thus Davies' Theorem says that
CH is equivalent to that every function $f:\R\times\R\to\R$ has
Davies representation.

The purpose of this paper is to prove the following descriptive
set-theoretic analogue of Davies' Theorem:

\thm{Theorem 1.}{Every $\Sigma^1_2$ function $f:\R\times\R\to\R$ has
a Davies representation
$$
f(x,y)=\sum_{n=0}^\infty g(x,n)h(y,n),
$$
where $g,h:\R\times\omega\to\R$ are $\Sigma^1_2$ functions and the
sum has only finitely many non-zero terms at each $(x,y)\in\R^2$, if
and only if all reals are constructible.}

\medskip

We will also show that it is not possible to find a Davies
representation of $f(x,y)=e^{xy}$ using Baire or Lebesgue measurable
functions $g$ and $h$. Note though that $e^{xy}$ {\it does} have a
representation as an infinite power series in $x$ and $y$. We will
give an example of a Borel (in fact, $\Delta^1_1$) function
$f:\R\times\R\to\R$ which does not admit a rectangular sum
representation as above with Baire or Lebesgue measurable $g$ and
$h$, even if we drop the pointwise finiteness condition of the sum,
and only ask that at each $(x,y)$ the sum converges pointwise.

\bigskip

{\it Organization}: In \S 2 below we show (Theorem 2) that if there
is a strongly $\Delta^1_n$ well-ordering of $\R$ then every
$\Sigma^1_n$ function $f:\R\times\R\to\R$ admits a representation
$$
f(x,y)=\sum_{n=0}^\infty g(x,n)h(y,n),
$$
with $\Sigma^1_n$ functions $g,h:\R\times\omega\to\R$, and where the
sum has only finitely many non-zero terms at each $(x,y)\in\R^2$.

In \S 3 we establish the converse to Theorem 2 in the case of
$\Sigma^1_2$ functions (Theorem 3). We also establish a converse in
the $\Sigma^1_3$ case, under the additional assumption that there is
a measurable cardinal. Finally, we establish the two facts regarding
representations using Baire and Lebesgue measurable functions
mentioned after Theorem 1 above.

\mysec{Inductive argument}

The necessary descriptive set-theoretic background for this paper
can be found in \cite{moschovakis} and \cite{mansfield}, in
particular the definitions of the (lightface) point-classes
$\Sigma^1_n$, $\Delta^1_n$ and $\Pi^1_n$. Here we recall the notions
for $\Delta^1_n$ well-orderings that are the most important to us.

Following \cite{harrington}, we say that a $\Delta^1_n$
well-ordering $\prec$ of $\R$ is {\it strongly} $\Delta^1_n$ if it
has length $\omega_1$ and the following (equivalent) statements hold
(c.f. \cite{moschovakis} chapter 5):

\begin{enumerate}
\item If $P\subseteq \R\times\R$ is $\Sigma^1_n$ then
$$
R(x,y)\iff (\forall z\prec y) P(x,z)
$$
is $\Sigma^1_n$.

\item The initial segment relation $\IS\subseteq \R\times
\R^{\leq\omega}$ defined by
$$
\IS(x,y)\iff (\forall z\prec x)(\exists n) y(n)=z\wedge (\forall
i,j) i=j\vee y(i)\neq y(j)
$$
is $\Sigma^1_n$.
\end{enumerate}

If all reals are constructible then there is a strongly $\Delta^1_2$
well-ordering of $\R$, see e.g. \cite{kanamori}.

It will often be necessary to work with recursively presented Polish
spaces other than $\R$, such as $\omega^\omega$ or $\R^{\leq\omega}$
(see below). Since all uncountable recursively presented Polish
spaces are isomorphic in the sense that there is a $\Delta^1_1$
bijection between them with a $\Delta^1_1$ inverse (see
\cite[3E.7]{moschovakis}), once we have a strongly $\Delta^1_n$
well-ordering of $\R$ we have a strongly $\Delta^1_n$ well-ordering
of all recursively presented Polish spaces. For convenience we will
use the same symbol, usually $\prec$, for such a well-ordering in
all the recursively presented spaces we consider. This minor
ambiguity poses no real danger.

We will say that a function $f:X\to Y$ from one recursively
presented Polish space $X$ to another, $Y$, is $\Sigma^1_n$
(respectively $\Pi^1_n$ and $\Delta^1_n$) if its graph is
$\Sigma^1_n$ (respectively $\Pi^1_n$ and $\Delta^1_n$). A function
$f:\R\times\R\to \R$ is said to have a {\it $\Sigma^1_n$ Davies
representation} if there are $\Sigma^1_n$ functions
$g,h:\R\times\omega\to\R$ such that
$$
\sum_{n=0}^\infty g(x,n)h(y,n)
$$
and the sum has only finitely many non-zero terms at each $(x,y)$.
The notions of $\Pi^1_n$ and $\Delta^1_n$ Davies representation are
defined similarly.

\thm{Theorem 2.}{If there is a strongly $\Delta^1_n$ well-ordering
of $\R$ then every $\Sigma^1_n$ function $f:\R\times\R\to\R$ has a
$\Sigma^1_n$ Davies representation. In particular, if all reals are
constructible then every $\Sigma^1_2$ function has a $\Sigma^1_2$
Davies representation.}

\medskip

To prove this, we will need to verify that Davies' proof, which uses
Zorn's Lemma, produces functions $g,h:\R\times\omega\to\R$ that are
$\Sigma^1_n$ and witness that $f$ has a $\Sigma^1_n$ Davies
representation. This in turn requires that we produce $\Sigma^1_n$
predicates (in the sense of \cite[p. 3]{mansfield} or \cite[p.
152--157]{kanamori}) that define $g$ and $h$. These predicates will
essentially be formulas defining $g$ and $h$ by transfinite
recursion as in the usual proof of the transfinite recursion
theorem, see e.g. \cite[p. 22, (2.6)]{jech}.

If $X$ is a set, we write $X^{\leq\omega}$ for the set of functions
$g: \alpha\to X$ for some $\alpha\in\omega+1$, and we set
$\lh(g)=|\dom(g)|$, the cardinality of $\dom(g)$. For $g\in
\R^{\leq\omega}$ we let
$$
\supp(g)=\{n\in\omega: n\in\dom(g)\wedge g(n)\neq 0\}.
$$

It is convenient for the proof to work relative to a fixed countable
sequence $x_n\in\mathscr P(\omega)$ of almost disjoint infinite
subsets of $\omega$. The sequence $(x_n)$ will be used to make sure
that certain almost disjoint families that are finite are not
maximal, because they will be constructed so that they are almost
disjoint from all $x_n$, $n\in\omega$. We will assume that the map
$n\mapsto x_n$ is recursive.

\medskip

{\it Definition.} The set $S\subseteq (\R^\omega)^{\leq\omega}\times
(\R^{\omega})^{\leq\omega}$ is defined by $(g,h)\in S$ if and only
if

\begin{enumerate}[(a)]
\item The sets $\supp(g(k))$, $\supp(h(m))$ and $x_n$ ($k\in\dom(g)$, $m\in\dom(h)$, $n\in\omega$) form an almost disjoint
sequence of sets.

\item For all $m\in\dom(g)$ there are infinitely many $k$ such that
$g(m)(k)=1$.

\item For all $n\in\dom(h)$ there are infinitely many $k$ such that
$h(n)(k)=1$.

\end{enumerate}

\bigskip

\noindent Note that $S$ is $\Delta^1_1$. We need the following Lemma
to encode the inductive step.

\mythm{Lemma.}{Suppose $f\in\R^{\leq\omega}$ is given and $(g,h)\in
S$ is such that $\lh(h)=\lh(f)$. Then there is
$\theta=\theta(f,g,h):\omega\to\R$ such that

\begin{enumerate}[\rm (1)]

\item For all $k\in\dom(f)$,
$$
f(k)=\sum_{l=0}^\infty \theta(l)h(k)(l),
$$
and the sum has only finitely many non-zero terms.
\item For all $n\in\dom(h)$, $\supp(\theta)\cap \supp(h(n))$ is
finite.
\item For all $n\in\dom(g)$, $\supp(\theta)\cap\supp(g(n))$ is
finite.
\item For all $n\in\omega$, $\supp(\theta)\cap x_n$ is finite.
\item For infinitely many $k$ we have $\theta(k)=1$.

\end{enumerate}

Moreover, $\theta$ may be found recursively in the given data. In
particular, there is a $\Delta^1_1$ function
$\theta:\R^{\leq\omega}\times S\to\R^\omega$ such that
$\theta(f,g,h)$ satisfies {\rm (1) -- (5)} for all
$(f,g,h)\in\R^{\leq\omega}\times S$.}

\begin{proof}

We define by induction on $k\in\omega$ an increasing sequence
$n_k\in\omega$ and $\theta\restrict n_k+1$ such that
\begin{enumerate}[(1')]

\item For all $m\in \dom(f)\cap (k+1)$,
$$
f(m)=\sum_{l=0}^{n_m} \theta(l)h(m)(l).
$$

\item For all $m\in \dom(h)\cap (k+1)$,
$\supp(\theta\restrict n_k+1)\cap\supp(h(m))\subseteq n_{m}+1$.

\item For all $m\in\dom(g)\cap (k+1)$,
$\supp(\theta\restrict n_k+1)\cap\supp(g(m))\subseteq n_{m}+1$.

\item For all $m\leq k$,
$\supp(\theta\restrict n_k+1)\cap x_m\subseteq n_{m}+1$.

\item $\theta(n_k)=1$.
\end{enumerate}
Assuming this can be done, $\theta$ will be defined on all of
$\omega$, since $n_k$ is increasing; By (1') and (2') it follows
that for $m\in\dom(f)$ we will have
$$
f(m)=\sum_{l=0}^{\infty} \theta(l)h(m)(l)
$$
and by (2') it holds that $\theta(l)h(m)(l)=0$ for $l>n_m$. Thus (1)
and (2) of the statement of the Lemma holds. Finally, (3'), (4') and
(5') ensures (3), (4) and (5).

To see that we can satisfy (1')--(5'), suppose $n_k$ and
$\theta\restrict n_k+1$ have been defined.

Case 1: $k+1\notin\dom(f)$. Then we let $p>n_k$ be the least number
greater than $n_k$ such that $p\notin\supp(g(m))$,
$p\notin\supp(h(m))$ and $p\notin x_m$, for $m\leq k$. The number
$p$ exists because of condition (a) in the definition of $S$. Define
$n_{k+1}=p$, and for $n_k<l< n_{k+1}$ let $\theta(l)=0$, and
$\theta(n_{k+1})=1$. Clearly conditions (1')--(5') are satisfied.

Case 2: $k+1\in\dom(f)$. Then let $p>n_k$ be the least number
greater than $n_k$ such that $p\notin\supp(g(m))$,
$p\notin\supp(h(m))$ and $p\notin x_m$, for $m\leq k$, and
$h(k+1)(p)=1$. The number $p$ exists because of conditions (a) and
(c) in the definition of $S$. We let $q>p$ be least such that
$q\notin\supp(g(m))$, $q\notin\supp(h(m))$ and $q\notin x_m$, for
$m\leq k+1$. Let $n_{k+1}=q$ and define for $n_k<l\leq n_{k+1}$,
$$
\theta(l)=\left\{ \begin{array}{ll}
1 & \text{if } l=q=n_{k+1},\\
f(k+1)-\sum_{m=0}^{p-1}\theta(m)h(k+1)(m)\ & \text{if } l=p,\\
0 & \textrm{otherwise.}
\end{array} \right.
$$
It is easy to see that (2')--(5') are satisfied. To see (1'), note
that
\begin{align*}
\sum_{l=0}^{n_{k+1}}\theta(l)h(k+1)(l)&=\sum_{l=0}^{p}\theta(l)h(k+1)(l)\\
&=f(k+1)-\sum_{m=0}^{p-1}\theta(m)h(k+1)(m)+\sum_{l=0}^{p-1}
\theta(l)h(k+1)(l)\\
&=f(k+1).
\end{align*}
This ends Case 2. It is clear from the construction that $\theta$ is
recursive in the given data $(f,g,h)$. Thus the map
$(f,g,h)\mapsto\theta(f,g,h)$ is in particular $\Delta^1_1$.
\end{proof}

\myheadpar{Davies' argument as an inductive
construction.\label{indconstr}} For the remainder of this section of
the paper, $\theta$ will be the function defined in Lemma 2.1. Using
this lemma one can now produce a Davies representation of
$f:\R\times\R\to\R$ by induction as follows: Assuming CH, fix a
well-ordering $\prec$ of $\R$ of order type $\omega_1$. Suppose
$g,h:\{y\in\R:y\prec x\}\times\omega\to\R$ have been defined such
that $\forall y,z\prec x$,
$$
f(y,z)=\sum_{n=0}^{\infty}g(y,n)h(z,n)
$$
and that further if $(w_m)$ is an enumeration of the initial segment
$\{y:y\prec x\}$ then the functions
$$
g_0(m)(n)=g(w_m,n)\text{ and } h_0(m)(n)=h(w_m,n)
$$
satisfy that $(g_0,h_0)\in S$. Then if we let $f_0(m)=f(w_m,x)$ and
define $g(x,n)=\theta(f_0,g_0,h_0)(n)$ then it is easy to check
using Lemma 2.1 that for $y\prec x$,
$$
f(x,y)=\sum_{n=0}^\infty g(x,n)h(y,n).
$$
If $(w'_m)$ enumerates $\{y:y\preceq x\}$ and we let
$f_1(m)=f(x,w'_m)$ and
$$
g_1(n)=\left\{ \begin{array}{ll}
\theta(f_0,g_0,h_0) & \text{if } w'_n=x\\
g_0(k) & \text{if } w'_n=w_k
\end{array} \right.
$$
then $(h_0,g_1)\in S$ and if we let
$h(x,n)=\theta(f_1,h_0,g_1))(n)$, it is again easy to check using
the previous Lemma that for all $y\preceq x$,
$$
f(y,x)=\sum_{n=0}^\infty g(y,n)h(x,n).
$$
Finally,
$$
g_1(m)(n)=g(w'_m,n)\text{ and } h_1(m)(n)=h(w'_m,n)
$$
satisfy $(g_1,h_1)\in S$, thus allowing the induction to continue.

\bigskip

Our task is now to verify that if $\prec$ is a strongly $\Delta^1_n$
well-ordering of $\R$, then the construction we have described may
be carried out in such a way that if $f:\R\times\R\to\R$ is
$\Sigma^1_n$, then the functions $g,h:\R\times\omega\to\R$ will be
$\Sigma^1_n$. This can be done since the strongly $\Delta^1_n$
well-ordering allows us to enumerate initial segments in a uniformly
$\Delta^1_n$ way. However, in order to be able to write down
$\Sigma^1_n$ definitions of $g$ and $h$ we need a lemma which says
that there is a $\Sigma^1_n$ function which can correctly compute
$g\restrict\{y:y\prec x\}\times\omega$ and $h\restrict\{y:y\prec
x\}\times\omega$ for every $x$.

Before stating that lemma we introduce various functions and
predicates. Fix a strongly $\Delta^1_n$ well-ordering $\prec$ of
$\R$ and let $\IS\subseteq \R\times\R^{\leq\omega}$ be the initial
segment relation as defined at the beginning of this section. Define
$\IS^*:\R\to\R^{\leq\omega}$ by
$$
\IS^*(x)=y\iff \IS(x,y)\wedge (\forall z\prec y)\neg\IS(x,z).
$$
Note that $\IS^*$ is $\Delta^1_n$. We also define a partial function
$\IS^\#:\R\times\R\to \omega$ by
$$
\IS^\#(x,y)=n\iff \IS^*(x)(n)=y.
$$
Note that the graph of $\IS^\#$ is a $\Delta^1_n$ subset of
$\R\times\R\times\omega$, and that if $y\prec x$ then $\IS^\#(x,y)$
computes the unique $n$ which $y$ corresponds to in the enumeration
of the initial segment of $x$ given by $\IS^*(x)$. Finally, we
define
$$
\Succ(x)=y\iff (\forall z\prec y) z=x\vee z\prec x.
$$

\mythm{Lemma.}{Let $f:\R\times\R\to \R$ be $\Sigma^1_n$ and suppose
there is a strongly $\Delta^1_n$ well-ordering $\prec$ of $\R$. Then
there is a unique $\Sigma^1_n$ function $F:\R\to(\R^\omega)^{\leq
\omega}\times (\R^\omega)^{\leq \omega}$ satisfying $F(x)=(G,H)$ if
and only if
\begin{enumerate}[\rm (1)]
\item $\lh(G)=\lh(H)=\lh(\IS^*(x))$ and $(G,H)\in S$.

\item If $z,z'\prec x$, $\IS^\#(x,z)=k$ and $\IS^\#(x,z')=k'$ then
$$
f(z,z')= \sum_{n=0}^\infty G(k)(n)H(k')(n).
$$

\item For all $y\prec x$, if we let
$w'=\IS^*(y)$, $w=\IS^*(x)$ and $f_0(k)=f(y,w'(k))$, and define for
$k\in\dom(w')$,
$$
G'(k)=G(l)\iff w'(k)=w(l)
$$
and
$$
H'(k)=H(l)\iff w'(k)=w(l)
$$
then $w(m)=y$ implies that
$$
G(m)=\theta(f_0,G',H')
$$

\item For all $y\prec x$, if we let $w'=\IS^*(y)$,
$w''=\IS^*(\Succ(y))$, $w=\IS^*(x)$ and $f_1(k)=f(w''(k),y)$, and
define for $k\in\dom(w'')$,
$$
G''(k)=G(l)\iff w''(k)=w(l)
$$
and for $k\in\dom(w')$,
$$
H'(k)=H(l)\iff w'(k)=w(l)
$$
then $w(m)=y$ implies that
$$
H(m)=\theta(f_1,H',G'')
$$
\end{enumerate}
}
\begin{proof}
Conditions (1)--(4) express {\it exactly} that for $y\prec x$, if we
let
$$
g(y,n)=G(\IS^\#(x,y))(n)
$$
and
$$
h(y,n)=H(\IS^\#(x,y))(n)
$$
then $g$ and $h$ are the functions we have constructed at stage $x$
in the inductive construction described in 2.2 above, if we at any
stage during the induction use the enumeration of the initial
segments given by the function $\IS^*$. Thus $F$ is unique and
defined for all $x$. Finally we note that the conditions (1)--(4)
can be expressed using $\Sigma^1_n$ predicates when $f$ is a
$\Sigma^1_n$ function. For instance, (3) may be replaced by
\begin{align*}
&(\forall y\prec x)(\exists w,w',f_0,G',H'\in
\R^{\leq\omega})(w'=\IS^*(y)\wedge w=\IS^*(x)\wedge\\
& \lh(f_0)=\lh(w')\wedge(\forall k\in\dom(w'))(
f_0(k)=f(y,w'(k))\wedge\\
& (\forall l\in\dom(G)) (G'(k)=G(l)\wedge H'(k)=H(l)\iff w'(k)=w(l)))\\
& \wedge (\forall m\in\dom(w)) (w(m)\neq y\vee
G(m)=\theta(f_0,G',H'))).
\end{align*}

 Thus (1)--(4) gives a $\Sigma^1_n$ definition
of the graph of $F$, and so the function $F$ is $\Sigma^1_n$.
\end{proof}

\begin{proof}[Proof of Theorem 2]
If $f:\R\times\R\to\R$ is $\Sigma^1_n$ and $\prec$ is a strongly
$\Delta^1_n$ well-ordering, let $F$ be as in the Lemma 2.3 and let
$F(x)=(G(x),H(x))$ for all $x$. Then
$$
g(x,n)=G(\Succ(x))(\IS^\#(\Succ(x),x))(n)
$$
and
$$
h(x,n)=H(\Succ(x))(\IS^\#(\Succ(x),x))(n).
$$
define $\Sigma^1_n$ functions that give us a Davies representation
of $f$.
\end{proof}

\remark{Remark.} If $f:\R\times\R\to\R$ is $\Delta^1_n$ then
conditions (1)--(4) define a $\Delta^1_n$ function $F$.
Consequently, the functions $g$ and $h$ produced in the proof of
Theorem 2 will be $\Delta^1_n$. Therefore we have:

\mythm{Corollary.}{If there is a strongly $\Delta^1_n$ well-ordering
of $\R$ then every $\Delta^1_n$ function $f:\R\times\R\to\R$ has a
$\Delta^1_n$ Davies representation.}

\bigskip

\mysec{A definable converse.}

We now show the following converse to Theorem 2 for $\Sigma^1_2$
functions:

\medskip

\thm{Theorem 3.}{If there are $\Sigma^1_2$ functions
$g,h:\R\times\omega\to\R$ such that
$$
e^{xy}=\sum_{n=0}^{\infty} g(x,n)h(y,n)
$$
with only finitely many non-zero terms at each $(x,y)$ then there is
a $\Sigma^1_2$ well-ordering of $\R$.}

\medskip

Since by Mansfield's Theorem (\cite{mansfield2}, see also
\cite[25.39]{jech}) the existence of a $\Sigma^1_2$ well-ordering or
$\R$ is equivalent to that all reals are constructible, Theorem 3
together with Theorem 2 proves Theorem 1 as stated in the
introduction. The proof requires several lemmata:

\mythm{Lemma.}{Let $b_0,\ldots, b_n\in\R$ be distinct reals and
$c_0,\ldots,c_n\in\R$. Then
$$
f(x)=\sum_{j=0}^n c_je^{xb_j}
$$
has $n+1$ distinct roots if and only if $c_0=\cdots=c_n=0$.}
\begin{proof}
By induction on $n$. If $f(x)$ has $n+1$ distinct roots then so does
$$
g(x)=e^{-b_0x}f(x).
$$
Using Rolle's Theorem from calculus it follows that $g'(x)$ has $n$
distinct roots, and so by the inductive hypothesis must be constant
zero. Thus $f(x)$ is constant zero.
\end{proof}

\mythm{Lemma.}{Let $a_0,\ldots, a_n$ and $b_0,\ldots, b_n$ be two
distinct sequences of real numbers. Then there are no functions
$g_l,h_l:\R\to\R$, $l<n$, such that
$$
e^{a_ib_j}=\sum_{l=0}^{n-1} g_l(a_i)h_l(b_j)
$$
for all $0\leq i,j\leq n$.}
\begin{proof}
If so then we have the matrix identity
$$
[e^{a_ib_j}]=\left[\begin{array} {ccc} g_0(a_0) &
\cdots & g_{n-1}(a_0)\\
\vdots & \ & \vdots\\
g_0(a_n) & \cdots & g_{n-1}(a_n)
\end{array}\right]\left[\begin{array} {ccc} h_0(b_0) &
\cdots & h_0(b_n)\\
\vdots & \ & \vdots\\
h_{n-1}(b_0) & \cdots & h_{n-1}(b_n)
\end{array}\right]
$$
and so $[e^{a_ib_j}]$ is a product of an $n+1\times n$ and an
$n\times n+1$ matrix. It follows that $\rank([e^{a_ib_j}])\leq n$,
which contradicts the previous Lemma.
\end{proof}

\mythm{Lemma.}{Assume $\Sigma^1_n$ uniformization holds and that
there are $\Sigma^1_n$ functions $g,h:\R\times\omega\to\R$ such that
$$
e^{xy}=\sum_{n=0}^{\infty} g(x,n)h(y,n)
$$
with only finitely many non-zero terms at each $(x,y)$. Suppose
there is an uncountable $\Sigma^1_n$ set $A\subseteq\R$ and a binary
$\Sigma^1_n$ relation $\prec$ on $\R$ such that $(A,\prec)$ is
well-ordered. Then there is a $\Sigma^1_n$ well-ordering of $\R$.}

\begin{proof}
Define
$$
N(x,y)=k\iff g(x,k)h(y,k)\neq 0\wedge (\forall l>k) g(x,l)h(y,l)=0.
$$
Clearly $N:\R\times\R\to\omega$ is $\Sigma^1_n$. Also define
$Q\subseteq\R\times\omega$ by
\begin{align*}
Q(x,n)\iff &(\exists a\in\R^\omega)(\forall i)(\forall j) (i=j\vee
a(i)\neq a(j))\wedge\\
&(\forall k) (a(k)\in A\wedge N(x,a(k))=n)\
\end{align*}
which is $\Sigma^1_n$. Let $Q^*:\R\to\omega$ be a $\Sigma^1_n$
uniformization of $Q$. Note that $Q^*$ is defined everywhere since
$A$ is uncountable.

Now define $R\subseteq\R\times[\R]^{<\omega}$, where
$[\R]^{<\omega}$ denotes the set of finite subsets\footnote{Formally
we let $[\R]^{<\omega}=\{s\in\R^{<\omega}: (\forall k<\lh(s)-1)
s(k)<s(k+1)\}$, where $<$ is the usual linear ordering of $\R$. Note
that for $s\in [\R]^{<\omega}$, the quantifiers $(\forall x\in s)$
and $(\exists x\in s)$ can be replaced by number quantifiers in
hierarchy calculations.} of $\R$, by
$$
R(x,s)\iff |s|=Q^*(x)+2\wedge (\forall y\in s) (y\in A\wedge
N(x,y)=Q^*(x)).
$$
Let $R^*:\R\to [\R]^{<\omega}$ be a $\Sigma^1_n$ uniformization of
$R$.

\claim{Claim.} $R^*$ is finite-to-1.

\begin{proof}
Suppose not. Then there is some $s=\{b_0,\ldots b_n\}$ such that
${R^*}^{-1}(s)$ is infinite. Pick $a_0,\ldots, a_n\in {R^*}^{-1}(s)$
distinct. Note that since $R^*(b_i)=s$ we have $Q^*(b_i)=|s|-2=n-1$.
Thus
$$
e^{a_ib_j}=\sum_{l=0}^{n-1} g(a_i,l)h(b_j,l),
$$
contradicting the previous Lemma.
\end{proof}

Let $\prec_{\lex}$ be the lexicographic order on $[A]^{<\omega}$
that we obtain from the well-ordering $\prec$ of $A$. Then we define
$<^*$ by
$$
x<^* y\iff R^*(x)\prec_{\lex} R^*(y)\vee (R^*(x)=R^*(y)\wedge x<y),
$$
where $<$ is the usual linear ordering of $\R$. Since $R^*$ is
finite-to-1, $<^*$ is a $\Sigma^1_n$ well-ordering of $\R$.
\end{proof}

\mythm{Lemma.}{There are no Baire or Lebesgue measurable
$g,h:\R\times\omega\to\R$ such that
$$
e^{xy}=\sum_{n=0}^\infty g(x,n)h(y,n)
$$
where the sum has finitely many non-zero terms at each $(x,y)$.}

\begin{proof}
Suppose there are Baire measurable $g,h:\R\times\omega\to\R$
representing $e^{xy}$ as above. Then
$$
N(x,y)=k\iff g(x,k)h(y,k)\neq 0\wedge (\forall l>k) g(x,l)h(y,l)=0.
$$
is also Baire measurable. It follows that there is some $n_0$ such
that
$$
A=\{(x,y)\in\R: N(x,y)=n_0\}
$$
is non-meagre and has the property of Baire. Thus we may find
$U,V\subseteq\R$ open and non-empty such that $A$ is comeagre in
$U\times V$. By Kuratowski-Ulam's Theorem it follows that
$$
\{x\in U: A_x \text{ is comeagre in } V\}
$$
is comeagre in $U$. Hence we may pick distinct elements $a_0,\ldots,
a_{n_0+1}\in U$ such that the section $A_{a_i}$ is comeagre in $V$
for all $i=0,\ldots, n_0+1$. But then we can find distinct elements
$$
b_0,\ldots, b_{n_0+1}\in\bigcap_{i=0}^{n_0+1} A_{a_i},
$$
which gives us that for $0\leq i,j\leq n_0+1$,
$$
e^{a_ib_j}=\sum_{n=0}^{n_0} g(a_i,n)h(b_j,n),
$$
contradicting Lemma 3.2.

\medskip

The proof of the Lebesgue measurable case is similar.
\end{proof}

\begin{proof}[Proof of Theorem 3.]
Suppose we have $\Sigma^1_2$ functions $g,h:\R\times\omega\to\R$
representing $e^{xy}$. By the previous Lemma, $g$ and $h$ cannot be
Baire measurable, and so $L\cap\R$ cannot be countable by
\cite[26.21]{jech}. But then we can apply Lemma 3.3 with $A=L\cap\R$
and $\prec$ the canonical $\Sigma^1_2$ well-ordering of $L\cap\R$ to
get a $\Sigma^1_2$ well-ordering of $\R$.
\end{proof}

\remark{Remark.} Assume $\Sigma^1_3$ uniformization. Suppose there
is a measurable cardinal and let $U$ be a normal ultrafilter
witnessing this. Then the tree representation for $\Sigma^1_3$ (see
\cite{kanamori} p. 201, also \cite[32.14]{jech}) and
\cite[15.10]{kanamori} gives us that if $\R\cap L[U]$ is countable
then all $\Sigma^1_3$ functions have the property of Baire. Since by
\cite[4.6]{silver} there is a $\Sigma^1_3$ well-ordering of $\R\cap
L[U]$, the proof above then gives us that if there is a $\Sigma^1_3$
Davies representation of $e^{xy}$ then there is a $\Sigma^1_3$
well-ordering of $\R$. In fact, we obtain the following stronger
result:

\mythm{Corollary.}{Assume $\Sigma^1_3$ uniformization. Suppose there
is a measurable cardinal and let $U$ be a normal ultrafilter
witnessing this. Then if there are $\Sigma^1_3$ functions
$g,h:\R\times\omega\to\R$ such that
$$
e^{xy}=\sum_{n=0}^{\infty} g(x,n)h(y,n)
$$
with only finitely many non-zero terms at each $(x,y)$ then
$\R=\R\cap L[U]$ and so there is a strongly $\Delta^1_3$
well-ordering of $\R$.}

\begin{proof}
By inspecting the proof of Lemma 3.3, there exists a finite-to-1
$\Sigma^1_3$ function $\theta:\R\to\R\cap L[U]$. Since the relation
$R\subseteq \R\times\N$ defined by
$$
R(y,n)\iff (\exists x_1,\ldots x_n) \bigwedge_{i=1}^n
\theta(x_i)=y\wedge\bigwedge_{i\neq j} x_i\neq x_j
$$
is $\Sigma^1_3$, it is absolute for transitive models containing
$U$. Suppose that there is $x_1\in \R\setminus L[U]$ and let
$y=\theta(x_1)$. If $n=|\theta^{-1}(y)\cap L[U]|$ then $R(y,n+1)$
holds in $V$. By absoluteness it holds in $L[U]$, contradicting that
$n=|\theta^{-1}(y)\cap L[U]|$. Thus $\R=\R\cap L[U]$ and by
\cite[5.2]{silver} there is a $\Delta^1_3$ wellordering of $\R$.
\end{proof}

\bigskip

In light of Theorem 3, it is natural to ask the following:

\remark{Question 1.} If there are $\Sigma^1_2$-functions $g_n,h_n$,
$n\in\omega$ such that
$$
e^{xy}=\sum_{n=0}^\infty g_n(x)h_n(y)
$$
with the sum having only finitely many non-zero terms at each
$(x,y)$, does the conclusion of Theorem 3 still hold? That is, is it
is necessary in Theorem 3 that $g_n,h_n$ are $\Sigma^1_2$ uniformly
in $n$?

\bigskip

In the paper \cite{shelah}, Shelah shows that the converse in
Davies' original Theorem does not remain true if we drop the
assumption that the sum must have at most finitely many non-zero
terms and only require the sum to converge pointwise. In a similar
vein we ask:

\remark{Question 2.} If we drop the finiteness condition, does
Theorem 1 still hold?

\bigskip

Shelah also shows in \cite{shelah} that if we add $\aleph_2$ Cohen
reals then there is a function $f:\R\times\R\to \R$ which does not
allow a representation
$$
f(x,y)=\sum_n^\infty g_n(x)h_n(y),
$$
even when we allow for the sum to have infinitely many non-zero
terms, requiring only that it converges pointwise. As a counterpoint
to that result, we point out the following:

\thm{Theorem 4.}{There is a Borel function $f:\R\times\R\to\R$ such
that for no $g_n,h_n:\R\to\R$ that are Baire measurable do we have
$$
f(x,y)=\sum_{n=0}^\infty g_n(x)h_n(y),
$$
for all $(x,y)\in\R^2$, where the sum converges point-wise, but may
have infinitely many non-zero terms. The same holds if we replace
Baire measurable by Lebesgue measurable.}

\begin{proof}
Let as usual $E_0$ denote the equivalence relation on $2^\omega$
defined by
$$
xE_0 y\iff (\exists N)(\forall n\geq N) x(n)=y(n).
$$
Let $\mathbf{1}_{E_0}$ be the characteristic function of $E_0$.
Suppose now that there are Baire measurable $g_n,h_n:\R\to\R$ such
that
$$
\mathbf{1}_{E_0}(x,y)=\sum_{n=0}^\infty g_n(x)h_n(y).
$$
Then we can find a dense $G_\delta$ set $A$ on which on which all
the functions $g_n$ and $h_n$ are continuous. But then for $x,y\in
A$ we have
$$
x E_0 y\iff (\forall k)(\exists N\geq k) \sum_{n=0}^N
g_n(x)h_n(y)>\frac 1 2.
$$
This gives us a $G_\delta$ definition of $E_0$ on $A$, and hence
$E_0$ must be a smooth equivalence relation on $A$ by
\cite{hakelou}, Corollary 1.2. But $E_0$ is not smooth on any
comeagre set, and we have a contradiction.

The proof of the Lebesgue measurable case is similar.
\end{proof}

\remark{Remark.} By \cite{judahshelah} (see also \cite[Exercise
26.2]{jech}), if there is a Cohen real (respectively random real)
over $L$ in $V$, then all $\Delta^1_2$ functions are Baire
measurable (respectively Lebesgue measurable). Thus it follows that
in this setting, $\mathbf 1_{E_0}$ cannot be represented as an
infinite pointwise convergent sum of rectangular $\Delta^1_2$
functions.

\bigskip

\begin{small}
{\sc\noindent Kurt G\"odel Research Center, University of Vienna\\
W\"ahringer Strasse 25, 1090 Vienna, Austria\\
{\it E-mail}: {\tt asger@logic.univie.ac.at}}

\bigskip

{\sc\noindent Department of Mathematics, University of Toronto\\
40 St. George Street, Room 6092, Toronto, Ontario, Canada\\
{\it E-mail}: {\tt weiss@math.utoronto.ca}}

\end{small}

\end{document}